\documentclass[oneside,11pt]{amsart}
\usepackage{float}
\usepackage{amscd,amssymb}
\usepackage[mathscr]{eucal}
\usepackage[all]{xy}
\usepackage{mathrsfs}
\usepackage{amsfonts}
\usepackage{amsmath}
\usepackage{amsthm}
\usepackage{latexsym}
\usepackage[all]{xy}

\usepackage{ifpdf}

\usepackage{graphicx}

\title{On the $\sf D$-affinity of quadrics in positive characteristic}

\author{Alexander Samokhin}

\address{Institute for Information Transmission Problems, Moscow, Russia}

\email{samohin@mccme.ru}

\thanks{This work was supported in part by the CNRS, the French Government
  Fellowship, and the RFFI award No. 02--01--22005.}

\newcommand{\Oo}{\mathcal O}
\newcommand{\Uu}{\mathcal U}

\newcommand{\bq}{\sf Q}
\newcommand{\Pp}{\mathbb P}

\newcommand{\Ff}{\mathcal F}

\newcommand{\Ll}{\mathcal L}

\newcommand{\D}{\mathcal D}

\newcommand*{\Dd}{\mathop{\mathrm D\kern0pt}\nolimits}
\newcommand*{\DD}{\mathop{\mathbb D\kern0pt}\nolimits}

\newcommand*{\Ext}{\mathop{\mathrm Ext}\nolimits}

\newtheorem{theorem}{Theorem}[section]
\newtheorem{corollary}{Corollary}[section]
\newtheorem{lemma}{Lemma}[section]

\newtheorem{remark}{Remark}[section]

\theoremstyle{definition}

\numberwithin{equation}{section}

\long\def\comment#1{}

\begin{document}

\begin{abstract}
In this paper we consider sheaves of differential operators on
quadrics of low dimension in positive characteristic. We prove a
vanishing theorem for the first term in the $p$-filtration of these
sheaves. This vanishing is a necessary condition for the $\sf D$-affinity of these
quadrics.   
\end{abstract}

\maketitle

\section{Introduction}
\label{}

Let $X$ be a smooth proper algebraic variety over an algebraically closed
field $k$ of arbitrary characteristic, ${\mathcal O}_X$ the structural
sheaf of $X$, and ${\mathcal D}_X$ the sheaf of differential operators
on $X$. Denote ${\mathcal M}({\mathcal D}_X)$ the category of (left)
${\mathcal D}_X$-modules that are coherent over $\D _X$. The
variety $X$ is said to be $\sf D$-affine if the two following
conditions hold: (i) for any ${\mathcal F}\in {\mathcal M}({\mathcal
  D}_X)$ one has ${\rm H}^{k}(X, \mathcal F) = 0$ for $k > 0$, and
(ii) the natural morphism $\D _X\otimes _{\Gamma (\D _X)}\Gamma
({\mathcal M})\rightarrow {\mathcal M}$ is surjective.

Beilinson and Bernstein proved (\cite{BB}) that homogeneous spaces of
semisimple algebraic groups are $\sf D$-affine if the base field $k$ has
characteristic zero. On the other hand, as was shown by 
Kashiwara et Lauritzen in \cite{KL}, homogeneous spaces over fields of positive characteristic 
can fail to be $\sf D$-affine. In this paper we
study quadrics of dimension less or equal to 4 in positive characteristic 
(these are homogeneous spaces of the orthogonal group). We prove 
a necessary condition for these quadrics to be $\sf D$-affine.
This theorem, which is the main result of the paper, is similar to
that of Andersen and Kaneda (\cite{AK}), where they treated the case of the flag
variety of the group in type ${\bf B}_2$.
 
Throughout we fix an algebraically closed field $k$ of characteristic
$p$. For a smooth variety $X$ over $k$ 
denote ${\sf F}\colon X\rightarrow X$ the absolute Frobenius morphism. Note that since $X$ is smooth, the
sheaf ${\sf F}_{\ast}{\mathcal O}_X$ is locally free. Recall that the sheaf of
differential operators ${\mathcal D}_X$ on $X$ admits the 
$p$-filtration defined as follows. For $r \geq 1$ denote ${\mathcal
  D}_r$ the endomorphism bundle ${\mathcal End}_{{\mathcal O}_X}({\sf
  F}^{r}_{\ast}{\mathcal O}_X)$ (here ${\sf F}^{r} = {\sf F}\circ
\dots \circ{\sf F}$ is the $r$-iteration of ${\sf F}$). Then
${\mathcal D}_X = \cup \ {\mathcal D}_r$. Recall that $X$ is said to be
Frobenius split if the sheaf ${\mathcal O}_X$ is a direct summand in
${\sf F}_{\ast}{\mathcal O}_X$. Homogeneous spaces of semisimple algebraic
groups are Frobenius split (\cite{MR}). For a Frobenius split variety $X$ and
for all $i\in \mathbb N$ the following property holds (Proposition, Sec. 1, \cite{AK}): 
\begin{equation}\label{eq:Lauritzen}
{\rm H}^{i}(X,{\mathcal D}_X) = 0 \Leftrightarrow {\rm H}^{i}(X,{\mathcal D}_r) = 0
\ \mbox{for any} \ r\in \mathbb N.
\end{equation}
\noindent Here is the main result of the paper 
(cf. the main theorem of \cite{AK}):

\vspace{0.1cm}

\begin{theorem}\label{th:mainth}
Let ${\sf Q}_n$ be a quadric of dimension $n\leq 4$. Then ${\rm H}^{i}({\sf Q}_n,{\mathcal
  D}_1) = 0$ for $i > 0$.
\end{theorem}

\section{Preliminaries}\label{sec:prel} 

Let $X$ be a smooth variety over $k$, and $\mbox{Coh}(X)$ the category
of coherent sheaves on $X$. The direct image functor
${\sf F}_{\ast}$ has a right adjoint functor ${\sf F}^{!}$ in
$\mbox{Coh}(X)$ (\cite{Har}). The duality theory for the finite flat
morphism ${\sf F}$ yields ({\it loc.cit.}):

\vspace{0.1cm}

\begin{lemma}\label{lem:rightadjoint}
The functor ${\sf F}^{!}$ is isomorphic to
\begin{equation}\label{eq:rightadjointequat}
{\sf F}^{!}(?) = {\sf F}^{\ast}(?)\otimes \omega _{X}^{1-p},
\end{equation}

\noindent where $\omega _X$ is the canonical invertible sheaf on $X$.
\end{lemma}

\vspace{0.1cm}

For any $i\geq 0$ one has an isomorphism,  the sheaf ${\sf F}_{\ast}{\mathcal O}_X$ being locally
free:

\begin{equation}
{\rm H}^{i}(X,{\mathcal End}_{{\mathcal O}_X}({\sf F}_{\ast}{\mathcal O}_X))
= \Ext ^{i}({\sf F}_{\ast}{\mathcal O}_X,{\sf F}_{\ast}{\mathcal O}_X).
\end{equation}
 
From Lemma \ref{lem:rightadjoint} we obtain:

\begin{eqnarray}
& \Ext ^{i}_{X}({\sf F}_{\ast}\Oo _{X},{\sf F}_{\ast}\Oo _{X}) = \Ext
^{i}({\Oo _{X}},{\sf F}^{!}{\sf F}_{\ast}(\Oo _{X})) =
{\rm H}^{i}(X,{\sf F}^{\ast}{\sf F}_{\ast}{\Oo _{X}}\otimes \omega _{X}^{1-p}).
\end{eqnarray}
Consider the fibered square:
\begin{center}
  \begin{center}
    \includegraphics[width=4cm]{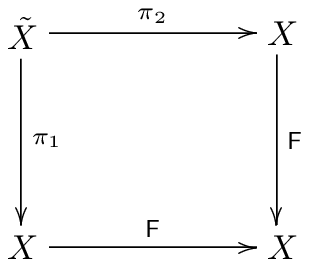}\label{fig:fig1}
  \end{center}
\end{center}

\vspace{0.1cm}

Let $i\colon \Delta \hookrightarrow X\times X$ be the diagonal
embedding, and $\tilde i$ the embedding $\tilde
X\hookrightarrow X\times X$ obtained from the above fibered
square. 

\vspace{0.1cm}

\begin{lemma}\label{lem:diagonal} 
One has an isomorphism of sheaves:

\begin{equation}
{\tilde i}_{\ast}\Oo _{\tilde X} = ({\sf F}\times {\sf F})^{\ast} (i_{\ast}\Oo _{\Delta}).
\end{equation}
\end{lemma}

Here ${\sf F}\times {\sf F}$ is the Frobenius morphism
on $X\times X$. The lemma is equivalent to saying that the fibered
product $\tilde X$ is isomorphic to the Frobenius neighbourhood of the
diagonal $\Delta \subset X\times X$.

\begin{proof}
Follows from the definition of fibered product.
\end{proof}

\vspace*{0.1cm}

\begin{lemma}\label{lem:isom}
There is an isomorphism of cohomology groups:
\begin{eqnarray}\label{eq:mainisom}
& {\rm H}^{i}(X,{\sf F}^{\ast}{\sf F}_{\ast}({\Oo _{X}})\otimes \omega _{X}^{1-p}) = 
{\rm H}^{i}(X\times X,({\sf F\times \sf F})^{\ast}(i_{\ast}{\Oo _{\Delta}})\otimes (\omega
  _{X}^{1-p}\boxtimes {\Oo _{X}})).
\end{eqnarray}
\end{lemma}

\noindent {\bf Proof.} Recall that the sign $\boxtimes$ in the right hand side of
(\ref{eq:mainisom}) denotes the external tensor product.
Applying the flat base change to the above fibered square, we get an
isomorphism of functors, the morphism $\sf F$ being flat:

\begin{equation}\label{eq:basechangefrob}
\sf F^{\ast}{\sf F}_{\ast} = {\pi _{1}}_{\ast}\pi _{2}^{\ast}.
\end{equation}

\noindent Note that all the functors ${\sf F}_{\ast}$, $\sf F^{\ast}$,
  ${\pi _{1}}_{\ast}$, and $\pi _{2}^{\ast}$ are exact, the morphism
  $\sf F$ being affine. The isomorphism (\ref{eq:basechangefrob}) implies an isomorphism of cohomology groups
\begin{equation}\label{eq:isom1}
{\rm H}^{i}(X,{\sf F}^{\ast}{\sf F}_{\ast}({\Oo _{X}})\otimes \omega _{X}^{1-p}) =
{\rm H}^{i}(X,{\pi _{1}}_{\ast}{\pi _{2}}^{\ast}({\Oo _{X}})\otimes \omega
_{X}^{1-p}).
\end{equation}

\noindent By the projection formula the right-hand side group in (\ref{eq:isom1}) is isomorphic to 
${\rm H}^{i}({\tilde X},{\pi _{2}}^{\ast}{\Oo _{X}}\otimes {\pi _{1}}^{\ast}{\omega
  _{X}^{1-p}})$. Let $p_1$ and $p_2$ be the projections of $X\times X$ onto
  the first and the second component respectively. One has
${\pi _1} = p_1\circ {\tilde i}$ and $\ {\pi _2} = p_2\circ
{\tilde i}$. Hence an isomorphism of sheaves
\begin{equation}\label{eq:isom2}
{\pi _{2}}^{\ast}{\Oo _{X}}\otimes {\pi _{1}}^{\ast}{\omega
  _{X}^{1-p}} = {\tilde i}^{\ast}(p_{2}^{\ast}{\Oo _{X}}\otimes
  p_{1}^{\ast}{\omega _{X}^{1-p}}) = {\tilde i}^{\ast}({\omega _{X}^{1-p}}
  \boxtimes {\Oo _{X}}).
\end{equation}

\noindent From (\ref{eq:isom2}) and the projection formula one obtains
\begin{equation}\label{eq:isom3}
 {\rm H}^{i}({\tilde X},{\pi
  _{2}}^{\ast}{\Oo _{X}}\otimes {\pi _{1}}^{\ast}{\omega
  _{X}^{1-p}}) = {\rm H}^{i}({\tilde X},{\tilde i}^{\ast}({\omega _{X}^{1-p}}
  \boxtimes {\Oo _{X}})) = {\rm H}^{i}(X\times X,{\tilde
  i}_{\ast}{\Oo _{\tilde X}}\otimes ({\omega _{X}^{1-p}}
  \boxtimes {\Oo _{X}})).
\end{equation}

\noindent Using Lemma \ref{lem:diagonal} we get the statement. \hfill$\Box$
 
\vspace{0.1cm}

\begin{corollary}\label{cor:yetanotherisom}

One has as well an isomorphism:

\begin{eqnarray}\label{eq:yetanotherisom}
{\rm H}^{i}(X,{\sf F}^{\ast}{\sf F}_{\ast}({\Oo _{X}})\otimes \omega
_{X}^{1-p}) = {\rm H}^{i}(X\times X,({\sf F}\times {\sf F})^{\ast}(i_{\ast}{\Oo _{\Delta}})\otimes ({\Oo  _{X}}\boxtimes
{\omega _{X}}^{1-p})).
\end{eqnarray}
\end{corollary}

Finally, we need a well-known lemma (e.g., {\rm SGA3}):

\vspace{0.1cm}

\begin{lemma}\label{lem:simpleusefullemma}
If a sheaf $\Ff$ on a variety $X$ is quasi-isomorphic to a
bounded complex $\Ff ^{\bullet}$ then ${\rm H}^{i}(X,\Ff) = 0$ provided that
${\rm H}^{p}(X,\Ff ^{q}) = 0$ for all $p + q = i$.
\end{lemma}

\section{Vanishing}

In this section we prove Theorem \ref{th:mainth}.
Let ${\sf Q}_n$ be a smooth quadric of dimension $n\leq 4$. 
Note that ${\bq}_1$ is isomorphic to $\Pp ^1$, and ${\bq}_2$ is isomorphic to
$\Pp ^1\times \Pp ^1$. For projective spaces in positive
characteristic the $\sf D$-affinity was proved in \cite{Haas} by
B. Haastert. We need, therefore, to consider the case $n=3$ and
$n=4$. We treat the case of three-dimensional quadrics, the 
four-dimensional case being similar. 
Let $\rm G$ be a simply connected simple algebraic group over $k$ of type
${\bf B}_2$, $\rm B$ a Borel subgroup of $\rm G$, and $\rm P \subset
\rm G$ a parabolic subgroup such that ${\rm G}/{\rm P}$ is isomorphic to a
quadric ${\sf Q}_3$. Denote $\pi \colon {\rm G}/{\rm B}\rightarrow {\rm G}/{\rm P}$
the projection. There exists a line bundle $\Ll$ over ${\rm G}/{\rm
  B}$ such that ${\rm R}^{0}\pi _{\ast}\Ll$ is a rank two vector bundle
over ${\sf Q}_3$, the spinor bundle. Denote $\Uu$ the dual bundle to
${\rm R}^{0}\pi _{\ast}\Ll$. There is a short exact sequence:

\begin{equation}\label{eq:tautseqon3dimquad}
0\rightarrow {\mathcal U}\rightarrow V\otimes \Oo _{{\bq}_3}\rightarrow {\mathcal
  U}^{\ast}\rightarrow 0,
\end{equation}

\noindent where $V$ is a symplectic $k$-vector space of dimension 4,
i.e. a space equipped with a non-degenerate skew form $\omega \in
\bigwedge ^{2}V^{\ast}$. 

\vspace{0.1cm}

\begin{lemma}\label{lem:quadricresofdiag}
Let ${\bq}_3$ be a smooth quadric of dimension $3$, and $i\colon
  \Delta \subset {\bq}_3\times {\bq}_3$ the diagonal embedding. Then the following
  complex is exact:

\begin{equation}\label{eq:resofdiag3dimquad}
0\rightarrow {\mathcal U}\boxtimes {\mathcal U}(-2)\rightarrow \Psi _2\boxtimes \Oo _{{\bq}_3}(-2)\rightarrow \Psi
_1\boxtimes \Oo _{{\bq}_3}(-1)\rightarrow \Oo _{{\bq}_3}\boxtimes \Oo
_{{\bq}_3}\rightarrow i_{\ast}\Oo _{\Delta}\rightarrow 0.
\end{equation}

\noindent  Here $\Psi _i$ for $i=1,2$ are some vector bundles on ${\bq }_3$, which
have right resolutions 
                                                    
\begin{equation}\label{eq:ampleresofkoszuldualonquadric}
0\rightarrow \Psi _{i}\rightarrow B_{i}\otimes _{k}\Oo
_{{\bq}_3}\rightarrow \dots \rightarrow B_{1}\otimes _{k}{\Oo
  _{{\bq}_3}}(i-1)\rightarrow {\Oo _{{\bq}_3}}(i)\rightarrow 0 ,
\end{equation} 


\noindent and $B_i$ are $k$-vector spaces.

\end{lemma}

\vspace{0.1cm}

\noindent {\bf Proof.} One can show, without recurrence to Kapranov's
theorem for quadrics (\cite{Kap}, Theorem 4.10), that the collection 
of bundles $\Uu (-2), \Oo _{{\bq}_3}(-2), \Oo _{{\bq}_3}(-1), \Oo _{{\bq}_3}$ is a
complete exceptional collection (\cite{Rud}) in $\Dd ^{b}({\bq}_3)$, 
the bounded derived category of coherent sheaves on ${\bq}_3$, 
and then use a purely categorical construction of a resolution of the diagonal ({\it
loc.cit.}) that can be performed over fields of arbitrary
characteristic. The bundles $\Uu , \Psi _2, \Psi _1, \Oo _{{\bq}_3}$
are terms of the so-called right dual collection. 
Otherwise, one can suitably modify Kapranov's argument
so that it will hold over fields of positive
characteristic. \hfill$\Box$

\vspace{0.1cm}

\begin{theorem}\label{th:mainthrevis}
One has  $\Ext ^{i}({\sf F}_{\ast}{\mathcal O}_{\bq _3},{\sf
  F}_{\ast}{\mathcal O}_{\bq _3}) = 0$ for $i > 0$.
\end{theorem}

\vspace{0.1cm}

\noindent {\bf Proof.} By Corollary \ref{cor:yetanotherisom} we need to show that 
${\rm H}^{i}({\bq}_3\times {\bq}_3,({\sf F\times \sf F})^{\ast}(i_{\ast}\Oo
_{\Delta})\otimes (\Oo  _{{\bq}_3}\boxtimes \omega _{{\bq}_3}^{1-p}))
= 0$ for $i>0$. Recall that $\omega _{{\bq}_3} = \Oo _{{\bq}_3}(-3)$. 
For a line bundle $\Ll$ one has ${\sf F}^{\ast}\Ll = \Ll
^p$. Taking the pull-back under $({\sf F}\times {\sf F})^{\ast}$
of the resolution (\ref{eq:resofdiag3dimquad}), we get 
a complex of coherent sheaves in degrees $-3,\dots ,0\colon$

\begin{equation}\label{eq:fpbresofdiag3dimquad}
0\rightarrow {\sf F}^{\ast}{\mathcal U}\boxtimes {\sf F}^{\ast}({\mathcal U}(-2))\rightarrow {\sf F}^{\ast}\Psi _2\boxtimes \Oo _{{\bq}_3}(-2p)\rightarrow {\sf F}^{\ast}\Psi
_1\boxtimes \Oo _{{\bq}_3}(-p)\rightarrow \Oo _{{\bq}_3}\boxtimes \Oo _{{\bq}_3}\rightarrow 0.
\end{equation}

\noindent Denote $C^{\bullet}$ the complex
(\ref{eq:fpbresofdiag3dimquad}), and let ${\tilde C}^{\bullet}$ be the tensor product of $C^{\bullet}$ with 
the invertible sheaf $\Oo  _{{\bq}_3}\boxtimes \omega
_{{\bq}_3}^{1-p}$. Then ${\tilde C}^{\bullet}$ is quasiisomorphic to
  the sheaf $({\sf F\times \sf F})^{\ast}(i_{\ast}\Oo _{\Delta})\otimes
  (\Oo _{{\bq}_3}\boxtimes \omega _{{\bq}_3}^{1-p})$. We thus have to
  compute the hypercohomology of ${\tilde C}^{\bullet}$. There is a distinguished
triangle in $\Dd ^{b}({\bq}_3)$:
\begin{equation}\label{eq:disttriangonquadric}
\dots \longrightarrow {\sf F}^{\ast}{\mathcal U}\boxtimes ({\sf
  F}^{\ast}{\mathcal U}^{\ast}\otimes \omega _{{\bq}_3})[3]\longrightarrow {\tilde C}^{\bullet}
\longrightarrow \sigma _{\geq -2}({\tilde
  C}^{\bullet})\stackrel{[1]}{\longrightarrow}\dots  \quad .
\end{equation}

\noindent Here $\sigma _{\geq -2}$ is the stupid truncation, and $[1]$
is a shift functor in $\Dd ^{b}({\bq}_3)$. Let us first look at the
  truncated complex $\sigma _{\geq -2}({\tilde C}^{\bullet})$, which is
  quasiisomorphic to 

\begin{equation}\label{eq:trunccomplex}
0\rightarrow {\sf F}^{\ast}\Psi _2\boxtimes \Oo _{{\bq}_3}(p-3) \rightarrow {\sf
  F}^{\ast}\Psi _1\boxtimes \Oo _{{\bq}_3}(2p-3) \rightarrow \Oo  _{{\bq}_3}\boxtimes
  \Oo _{{\bq}_3}(3p-3)\rightarrow 0.
\end{equation}

\noindent Let us show that ${\rm H}^{k}({\bq}_3,{\sf F}^{\ast}\Psi _i) = 0$ for
$k>i$ and $i=1,2$. Indeed, the sheaves $\Psi _1$ and $\Psi _2$ have
resolutions as in (\ref{eq:ampleresofkoszuldualonquadric}). By the
Kempf vanishing theorem (\cite{Ke}), effective line bundles on
homogeneous spaces have no higher cohomology. The terms of resolutions
(\ref{eq:ampleresofkoszuldualonquadric}) for $n=3$ and $i=1,2$ consist of direct sums of
effective line bundles and of direct sums of the sheaf $\Oo _{{\bq}_3}$.
Positivity of a line bundle is preserved under the Frobenius
pullback, hence the terms of the resolutions of sheaves ${\sf F}^{\ast}\Psi
_i$ have only zero cohomology. Applying Lemma
\ref{lem:simpleusefullemma}, we obtain the above vanishing.
Further, line bundles that occur in the second argument of the
terms of the complex (\ref{eq:trunccomplex}) are effective for $p>3$. 
For $p=2$ and $p=3$ the line bundle occurring in the leftmost term of
(\ref{eq:trunccomplex}) is isomorphic to
$\Oo _{{\bq}_3}$ and $\Oo _{{\bq}_3}(-1)$, respectively. For all $p$
these line bundles have no higher cohomology. Using again Lemma
\ref{lem:simpleusefullemma}, we get ${\mathbb H}^{i}(\sigma _{\geq
  -2}({\tilde C}^{\bullet})) = 0$ for $i>0$. 

Consider now the bundle ${\sf F}^{\ast}{\mathcal U}\boxtimes ({\sf F}^{\ast}{\mathcal
  U}^{\ast}\otimes \omega _{{\bq}_3})$. The proof of Theorem \ref{th:mainthrevis} will be
completed if we show that ${\rm H}^{i}({\bq}_3\times {\bq}_3,{\sf F}^{\ast}{\mathcal U}\boxtimes ({\sf F}^{\ast}{\mathcal
  U}^{\ast}\otimes \omega _{{\bq}_3})) = 0$ for $i\neq 3$. Indeed, 
taking cohomology of the triangle (\ref{eq:disttriangonquadric}), we see that ${\mathbb
  H}^{i}({\tilde C}^{\bullet}) = 0$ for $i>0$, q.e.d.  \hfill$\Box$.

\vspace{0.1cm}

\begin{lemma}\label{lem:mainlemma}
One has ${\rm H}^{i}({\bq}_3\times {\bq}_3,{\sf F}^{\ast}{\mathcal U}\boxtimes ({\sf F}^{\ast}{\mathcal
  U}^{\ast}\otimes \omega _{{\bq}_3})) = 0$ for $i\neq 3$. 
\end{lemma}

\vspace{0.1cm}

\noindent {\bf Proof.} By the Serre duality one has ${\rm H}^{k}({\bq}_3,{\sf F}^{\ast}{\mathcal U})$ = ${\rm H}^{3-k}({\bq}_3,{\sf F}^{\ast}{\mathcal U}^{\ast}\otimes \omega _{{\bq}_3})$,
for $0\leq k\leq 3$. Using the K\"unneth formula, we see that it is
sufficient to show that ${\rm H}^{k}({\bq}_3,{\sf F}^{\ast}{\mathcal U})$ is
  non-zero for just one value of $k$. In fact, ${\rm H}^{k}({\bq}_3,{\sf
      F}^{\ast}{\mathcal U}) = 0$ if $k\neq 2$. To prove this, apply
    the Frobenius pullback ${\sf F}^{\ast}$ to
    (\ref{eq:tautseqon3dimquad}). We get

\begin{equation}\label{eq:pulbcktautseqon3dimquad}
0\rightarrow {\sf F}^{\ast}{\mathcal U}\rightarrow {\sf F}^{\ast}V\otimes \Oo
_{{\bq}_3}\rightarrow {\sf F}^{\ast}{\mathcal U}^{\ast}\rightarrow 0
\end{equation}

\noindent We need a particular case of the following theorem due to
Carter and Lusztig (\cite{Ar}, Theorem 6.2):

\vspace{0.1cm}

\begin{theorem}\label{th:Lusztigtheor}
Let $k$ be a field of characteristic $p$ and let $E$ be a vector space
of dimension $r$ over $k$. Then there exists an exact sequence of $GL_{r}(k)$-modules 
\begin{equation}\label{eq:Lusztigres}
0\rightarrow {\sf F}^{\ast}E\rightarrow {\sf S}^{p}(E)\rightarrow \Sigma
^{(p-1,1)}(E)\rightarrow \dots \rightarrow \Sigma
^{\lambda}(E)\rightarrow 0
\end{equation}

\noindent where $\lambda = (p - \mbox{min} \
(p-1,r-1),1,1,\dots)$. Here $\Sigma ^{\lambda}$ are Schur functors. In
particular, if a partition $\lambda$ is equal to $(p,0,\dots)$ then
$\Sigma ^{\lambda}$ is equal to $p$-th symmetric power functor ${\sf
  S}^{p}$. This construction globalizes to produce a resolution of the
Frobenius pull-back of a vector bundle (\cite{Ar}). 
\end{theorem}

\vspace{0.1cm}

Theorem \ref{th:Lusztigtheor}, applied to the bundle $\Uu ^{\ast}$, furnishes a
short exact sequence, the bundle $\Uu ^{\ast}$ having rank 2:

\begin{equation}\label{eq:Lusztigreson3dimquad}
0\rightarrow {\sf F}^{\ast}{\mathcal U}^{\ast}\rightarrow {\sf S}^{p}{\mathcal
  U}^{\ast}\rightarrow {\sf S}^{p-2}{\mathcal U}^{\ast}\otimes \Oo
  _{{\bq}_3}(1)\rightarrow 0
\end{equation}

\noindent The vector bundles ${\sf S}^{p}{\mathcal U}^{\ast}$ and ${\sf
  S}^{p-2}{\mathcal U}^{\ast}\otimes \Oo _{{\bq}_3}(1)$ are
  pushforwards onto ${\bq}_3$ of effective line bundles $\Ll _{\chi
  _1}$ and $\Ll _{\chi _1}$ over ${\rm G}/{\rm B}$ that correspond to the weights
  $\chi _1 = (p,0)$ and $\chi _2 = (p-1,1)$, respectively. Moreover, the
  restrictions of $\Ll _{\chi _1}$ and $\Ll _{\chi _1}$ to the fibers of
  $\pi$ have positive degrees. Considering the Leray spectral sequence
  for the morphism $\pi$ and using the Kempf vanishing theorem, we
  obtain that the bundles ${\sf S}^{p}{\mathcal U}^{\ast}$ and ${\sf
  S}^{p-2}{\mathcal U}^{\ast}\otimes \Oo _{{\bq}_3}(1)$ both have no
  higher cohomology. By Lemma \ref{lem:simpleusefullemma} one then has 
  ${\rm H}^{i}({\bq}_3,{\sf F}^{\ast}{\mathcal U}^{\ast}) =
  0$ for $i>1$. Moreover, there is an isomorphism ${\rm H}^{0}({\bq}_3,{\sf
  F}^{\ast}{\mathcal U}^{\ast}) = V$. Considering the long exact cohomology
  sequence associated to (\ref{eq:pulbcktautseqon3dimquad}), we get the
  statement of Lemma \ref{lem:mainlemma}. Hence, it remains to prove the
  isomorphism  ${\rm H}^{0}({\bq}_3,{\sf F}^{\ast}{\mathcal U}^{\ast}) =
  V$. Indeed, taking cohomology of (\ref{eq:Lusztigreson3dimquad}), we
  get a short exact sequence:

\begin{equation}\label{eq:cohseqof364}
0\rightarrow {\rm H}^{0}({\bq}_3,{\sf F}^{\ast}{\mathcal U}^{\ast})\rightarrow {\rm H}^{0}({\bq}_3,{\sf S}^{p}{\mathcal
  U}^{\ast})\rightarrow {\rm H}^{0}({\bq}_3,{\sf S}^{p-2}{\mathcal U}^{\ast}\otimes \Oo
  _{{\bq}_3}(1))
\end{equation}

\noindent Using again the Kempf vanishing theorem, we obtain
  ${\rm H}^{0}({\bq}_3$,${\sf S}^{p}{\mathcal
  U}^{\ast}) = {\sf S}^{p}V^{\ast}$ and ${\rm H}^{0}({\bq}_3,{\sf
  S}^{p-2}{\mathcal U}^{\ast}\otimes \Oo
  _{{\bq}_3}(1)) = \Sigma ^{(p-1,1)}V^{\ast}$. On the other hand, the
  resolution (\ref{eq:Lusztigres}), applied to the vector space $V^{\ast}$,
  furnishes a short exact sequence:

\begin{equation}\label{eq:Lusztigresrk2}
0\rightarrow {\sf F}^{\ast}V^{\ast}\rightarrow {\sf S}^{p}V^{\ast}\rightarrow \Sigma
^{(p-1,1)}V^{\ast}
\end{equation}
 
\noindent Comparing (\ref{eq:cohseqof364}) and (\ref{eq:Lusztigresrk2}),
and taking into account the above isomorphisms, we get the isomorphism 
${\rm H}^{0}({\bq}_3,{\sf F}^{\ast}{\mathcal U}^{\ast}) = V$. Finally, one
has an isomorphism $V = V^{\ast}$ since $V$ is symplectic. This
implies Lemma \ref{lem:mainlemma}.  \hfill$\Box$

\vspace{0.1cm}

The proof of Theorem \ref{th:mainth} in the case of four-dimensional
quadrics is essentially the same, the only
difference being that there are two spinor bundles in this case. 

\begin{remark}
{\rm A more general result, as well as applications to derived categories
of coherent sheaves, are discussed in a forthcoming paper (\cite{Sam}).}
\end{remark}



\section*{Acknowledgements}
I am happy to thank R. Bezrukavnikov and A. Kuznetsov for extremely
valuable discussions. This paper was written a year ago during the
author's visit to the University Paris 13. I would like to express 
my gratitude to this institution for its hospitality.

\end{document}